\documentclass[12pt,a4paper,eqno]{amsart}
\usepackage{graphicx,epsfig}
\usepackage{amssymb}

\def\qed{\hfill{ $\clubsuit$ }}

\newtheorem{definition}{\bf Definition}[section]

\newenvironment{namelist}[1]{%
\begin{list}{}
     {
      
      \settowidth{\labelwidth}{#1}
      \setlength{\leftmargin}{1.1\labelwidth}
               }
      }{%
\end{list}}

\newtheorem{thm}{Theorem}[section]

\newtheorem{lem}[thm]{Lemma}
\newtheorem{prop}[thm]{Proposition}
\theoremstyle{definition}

\theoremstyle{remark}

\begin{document}

\title[Kendall random walk]%
      {Kendall random walk, Williamson transform and the corresponding Wiener-Hopf factorization}
\author{B.H. Jasiulis-Go{\l}dyn$^1$ and J.K. Misiewicz $^2$}
\thanks{$^1$ Institute of Mathematics, University of Wroc{\l}aw, pl. Grunwaldzki 2/4, 50-384 Wroc{\l}aw, Poland, e-mail: jasiulis@math.uni.wroc.pl \\
$^2$ Faculty of Mathematics and Information Science, Warsaw
University of Technology, ul. Koszykowa 75, 00-662 Warsaw, Poland, e-mail:
j.misiewicz@mini.pw.edu.pl 
           }
\maketitle

\begin{abstract}
The paper gives some properties of hitting times and an analogue of the Wiener-Hopf factorization for the Kendall random walk. We show also that the Williamson transform is the best tool for problems connected with the Kendall generalized convolution.
\end{abstract}

\section{Introduction}
In this paper we study positive and negative excursions for the Kendall random walk $\{X_n \colon n \in \mathbb{N}_0\}$ which can be defined by the following recurrence construction:
\begin{definition} The stochastic process $\{X_n \colon n \in \mathbb{N}_0\}$ is a discrete time Kendall random walk  with parameter $\alpha>0$ and step distribution $\nu$ if there exist 
\begin{namelist}{ll}
\item[\bf 1.] $(Y_k)$ i.i.d. random variables with distribution $\nu$,
\item[\bf 2.] $(\xi_k)$ i.i.d. random variables with uniform distribution on $[0,1]$,
\item[\bf 3.] $(\theta_k)$ i.i.d. random variables with symmetric Pareto distribution with density $\widetilde{\pi}_{2\alpha}\, (dy) = \alpha |y|^{-2\alpha-1} \pmb{1}_{[1,\infty)}(|y|)\, dy$,
\item[\bf 4.]  sequences $(Y_k)$, $(\xi_k)$ and $(\theta_k)$ are independent
\end{namelist}
such that
$$
X_0 = 1, \quad X_1 = Y_1, \quad  X_{n+1} = M_{n+1} r_{n+1} \left[ \mathbf{I}(\xi_n < \varrho_{n+1}) + \theta_{n+1} \mathbf{I}(\xi_n > \varrho_{n+1})\right],
$$
where $\theta_{n+1}$ and $M_{n+1}$ are independent,
$$
M_{n+1} = \max\{ |X_n|,|Y_{n+1}|\}, \quad m_{n+1} = \min\{ |X_n|,|Y_{n+1}|\}, \quad \varrho_{n+1} = \frac{m_{n+1}^{\alpha}}{M_{n+1}^{\alpha}}
$$ 
and
$$
r_{n+1} = \left\{sgn(u) : \max \left\{|X_n|,|Y_{n+1}| \right\} = |u| \right\}.
$$
\end{definition}
The Kendall random walk is an example of discrete time L\'{e}vy random walk under generalized convolution studied in \cite{BJMR}. Some basic properties of this particular  process are described in  \cite{KendallWalk}. 

In the second section we present the definition and main properties of the Williamson transform together with unexpectedly simple inverse formula. We give also the definition of the Kendall convolution (\cite{misjas1}) which is an example of generalized convolution in the sense defined by Urbanik \cite{Urbanik93, Urbanik64}. This section is written only for the completeness of the paper, since all presented facts are known, but rather difficult to find in the literature.

In the third section we recall another definition of Kendall random walk following the multidimensional distributions approach given in \cite{BJMR}. By this construction it is clear that the Kendall random walk is a Markov process. Then we study hitting times, positive and negative excursions for this process. The main result of the paper is the analog of the Wiener-Hopf factorization.

For simplicity we use notation $T_a$ for the rescaling operator (dilatation) defined by $(T_a\lambda)(A) = \lambda ({A/a})$ for every Borel set $A$ when $a \neq 0$, and $T_0 \lambda = \delta_0$. In the whole paper parameter $\alpha > 0$ is fixed.

\section{Williamson transform and  Kendall convolution}
Let us start from reminding the Williamson transform (see \cite{Williamson}) and its basic properties. In this paper we apply it to $\mathcal{P}_s$-symmetric probability measures on $\mathbb{R}$, but equivalently it can be considered on the set $\mathcal{P}_+$- measures on $[0,\infty)$.
\begin{definition}
By Williamson transform we understand operation $\nu \rightarrow \widehat{\nu}$ given by
$$
\widehat{\nu}(t) = \int_{\mathbb{R}} \left( 1 - |xt|^{\alpha} \right)_+ \nu(dx), \quad \nu \in \mathcal{P}_s,
$$
where $a_+ = a$ if $a\geqslant 0$ and $a_+=0$ otherwise.
\end{definition}

For the convenience we use the following notation:
$$
\Psi(t) = \left( 1 - |t|^{\alpha} \right)_+, \quad G(t) = \widehat{\nu}(1/t).
$$
The next lemma is almost evident and well known. It contains the inverse of Williamson transform which is surprisingly simple.
\begin{lem}
The correspondence between a measure $\nu \in \mathcal{P}$ and its Williamson transform is $1-1$. Moreover, denoting by $F$ the cumulative distribution function of $\nu$, $\nu(\{0\}) = 0$, we have
$$
F(t) = \left\{ \begin{array}{lcl}
 \frac{1}{2\alpha} \left[ \alpha (G(t) + 1) + t G'(t) \right] & \hbox{if} & t>0; \\
 1 - F(-t) & \hbox{if} & t<0.
 \end{array} \right.
$$
except for the countable many $t \in \mathbb{R}$.
\end{lem}

\noindent {\bf Proof.} Of course the assumption $\nu(\{0\}) = 0$ is only technical simplification. Notice first, that for $t>0$ integrating by parts we obtain
\begin{eqnarray*}
G(t) & = & \int_{\mathbb{R}} \left( 1 - |x/t|^{\alpha} \right)_+ dF(x) \\
 & = & 2 \left[ \bigl( 1 - x^{\alpha} t^{-\alpha}\bigr) F(x)\big|_0^t + \alpha t^{-\alpha} \int_0^t x^{\alpha - 1} F(x) dx \right] \\
 & = & -1 + 2\alpha t^{-\alpha} \int_0^t  x^{\alpha - 1} F(x) dx.
\end{eqnarray*}
Since
$$
\lim_{h\rightarrow 0^+} \frac{1}{h} \left[ \int_0^{t+h} x^{\alpha - 1} F(x) dx - \int_0^t x^{\alpha - 1} F(x) dx \right] = t^{\alpha-1} F(t+),
$$
and
$$
\lim_{h\rightarrow 0^-} \frac{1}{h} \left[ \int_0^{t+h} x^{\alpha - 1} F(x) dx - \int_0^{t} x^{\alpha - 1} F(x) dx \right] = t^{\alpha-1} F(t-)
$$
then, except for the points of jumps for $F$ we can differentiate the equality
$$
2\alpha  \int_0^t  x^{\alpha - 1} F(x) dx = t^{\alpha}\left( G(t) + 1 \right)
$$
and obtain
$$
2\alpha  t^{\alpha - 1} F(t) = \alpha t^{\alpha-1}\left( G(t) + 1 \right) + t^{\alpha} G'(t). \eqno{\clubsuit}
$$
In the following we use the notation: $\widetilde{\delta}_x = \frac{1}{2} \delta_x + \frac{1}{2} \delta_{-x}$ and $\widetilde{\pi}_{2\alpha}\, (dy) = \alpha |y|^{-2\alpha-1} \pmb{1}_{[1,\infty)}(|y|)\, dy$ is the symmetric Pareto distribution.
\begin{definition}
By the Kendall convolution we understand operation $\vartriangle_{\alpha} \colon \mathcal{P}_s^2 \rightarrow \mathcal{P}_s$ defined for discrete measures by
$$
\widetilde{\delta}_x \vartriangle_{\alpha} \widetilde{\delta}_y := T_M \left( \varrho^{\alpha} \widetilde{\pi}_{2\alpha} + (1- \varrho^{\alpha}) \widetilde{\delta}_1 \right)
$$
where $M = \max\{|x|, |y|\}$, $m = \min\{|x|, |y|\}$, $\varrho = {m/M}$. The extension of $\vartriangle_{\alpha}$ to the whole $\mathcal{P}_s$ is given by
$$
\nu_1 \vartriangle_{\alpha} \nu_2 (A)  = \int_{\mathbb{R}^2} \widetilde{\delta}_x \vartriangle_{\alpha} \widetilde{\delta}_y (A)\, \nu_1(dx) \nu_2(dy).
$$
\end{definition}
It is easy to see that the operation $\vartriangle_{\alpha}$ is symmetric, associative, commutative and
\begin{itemize}
\item $\nu \vartriangle_{\alpha} \delta_0 = \nu$ for each $\nu \in \mathcal{P}_s$;
\item $\left( p \nu_1 + (1-p)\nu_2 \right) \vartriangle_{\alpha} \nu = p
\left( \nu_1 \vartriangle_{\alpha} \nu \right) + (1-p) \left( \nu_2 \vartriangle_{\alpha}
\nu \right)$ for each $p \in [0,1]$ and each $\nu, \nu_1, \nu_2 \in \mathcal{P}_s$.
\item if $\lambda_n \rightarrow \lambda$ and $\nu_n \rightarrow \nu$
then $(\lambda_n \vartriangle_{\alpha} \nu_n) \rightarrow (\lambda \vartriangle_{\alpha} \nu)$, where $\rightarrow$ denotes the weak convergence;
\item $T_a \bigl( \nu_1 \vartriangle_{\alpha}\nu_2 \bigr)
= \bigl(T_a \nu_1\bigr) \vartriangle_{\alpha} \bigl(
T_a \nu_2\bigr)$ for each $\nu_1, \nu_2 \in \mathcal{P}_s$.
\end{itemize}
The next proposition, which is only a reformulation of a known property of Williamson transform, shows that this transform is playing the same role for Kendall convolution as the Fourier transform for classical convolution.
\begin{prop}\label{prop:1}
Let $\nu_1, \nu_2 \in \mathcal{P}_s$ be probability measures with the Williamson transforms $\widehat{\nu_1}, \widehat{\nu_2}$. Then
$$
\int_{\mathbb{R}} \Psi(xt) \bigl(\nu_1 \vartriangle_{\alpha} \nu_2 \bigr)(dx) = \widehat{\nu_1}(t) \widehat{\nu_2}(t).
$$
\end{prop}

\noindent {\bf Proof.}
Notice first that for $M = \max\{|x|, |y|\}$, $m = \min\{|x|, |y|\}$, $\varrho = {m/M}$ we have
\begin{eqnarray*}
\lefteqn{\int_{\mathbb{R}} \Psi(ut) \bigl(\widetilde{\delta}_x \vartriangle_{\alpha} \widetilde{\delta}_y \bigr)(du) = \int_{\mathbb{R}} \Psi(Mut) \bigl( \widetilde{\delta}_{\varrho} \vartriangle_{\alpha} \widetilde{\delta}_1\bigr) (du)} \\
 & = & \varrho^{\alpha} 2 \int_0^{\infty}\!\!\Psi(Mut) \frac{\alpha\, du}{ |u|^{2\alpha+1}} + (1- \varrho^{\alpha}) \Psi(Mt) \mathbf{1}_{[-1,1]}(Mt) \\
 & = & \varrho^{\alpha} \left( 1 - |Mt|^{\alpha} \right)^2 \mathbf{1}_{[-1,1]}(Mt) + (1- \varrho^{\alpha}) \left( 1 - |Mt|^{\alpha} \right) \mathbf{1}_{[-1,1]}(Mt) \\[1mm]
 & = & \left( 1 - |Mt|^{\alpha} \right) \left( 1 - |mt|^{\alpha} \right) \mathbf{1}_{[-1,1]}(Mt) = \Psi(xt) \Psi(yt)
\end{eqnarray*}
Now, by the definition of Kendall convolution, we see that
\begin{eqnarray*}
\lefteqn{\hspace{-15mm} \int_{\mathbb{R}}\! \Psi(xt) \!\bigl(\nu_1 \!\vartriangle_{\alpha}\! \nu_2 \bigr)(dx)= \int_{\mathbb{R}^3} \Psi(ut) \bigl(\widetilde{\delta}_x \vartriangle_{\alpha} \widetilde{\delta}_y \bigr)(du) \, \nu_1(dx) \nu_2 (dy)} \\
& = & \int_{\mathbb{R}^2}  \Psi(xt) \Psi(yt) \nu_1(dx) \nu_2 (dy) = \widehat{\nu_1}(t) \widehat{\nu_2}(t),
\end{eqnarray*}
which was to be shown. \qed
\begin{prop}\label{prop:2}
There exists a sequence of positive numbers $(c_n)$ such that $%
T_{c_n} \widetilde{\delta}_1^{\vartriangle_{\alpha} n}$ converges weakly to a measure $\sigma \neq \delta_0$ (here $\nu^{ \vartriangle_{\alpha} n} = \nu \vartriangle_{\alpha} \dots \vartriangle_{\alpha} \nu$ denotes the Kendall convolution of $n$ measures $\nu$).
\end{prop}

\noindent {\bf Proof.} Since the Williamson transform of $\widetilde{\delta}_1$ is equal $\Psi(t)$ then putting $c_n = n^{-{1/{\alpha}}}$ we obtain
\begin{eqnarray*}
\lefteqn{\hspace{-20mm} \int_{\mathbb{R}}\! \Psi(xt)  T_{c_n} \widetilde{\delta}_1^{\vartriangle_{\alpha} n} (dx) = \int_{\mathbb{R}}\! \Psi(c_n x t) \widetilde{\delta}_1^{\vartriangle_{\alpha} n} (dx)} \\
& = & \left( 1 - |c_n t|^{\alpha} \right)_+^n \longrightarrow e^{-|t|^{\alpha}} \quad \hbox{if} \quad n \rightarrow \infty.
\end{eqnarray*}
The function $e^{-|t|^{\alpha}}$ is symmetric and  monotonically decreasing as a function of $t^{\alpha}$ on the positive half line, thus there exists a measure $\sigma \in \mathcal{P}_s$ such that
$$
e^{-|t|^{\alpha}} = \int_{\mathbb{R}}\! \left( 1 - |xt|^{\alpha} \right)_+  \sigma(dx).
$$
According to Lemma 2.1 we can calculate the cumulative distribution function $F$ of the measure $\sigma$ obtaining $F(0) = \frac{1}{2}$ and
$$
F(t) = \left\{ \begin{array}{lcl}
   \frac{1}{2} \bigl( 1 + t^{-\alpha} + e^{t^{-\alpha}} \bigr) e^{- t^{-\alpha}} & \hbox{if} & t> 0; \\
   1 - F(-t) & \hbox{if} & t< 0. \end{array} \right. \eqno{\clubsuit}
$$
{\bf Remark.} Notice that properties of Kendall convolution together with the one described in the last proposition show that the Kendall convolution is an example of  generalized convolutions in the sense defined by K. Urbanik. More about generalized convolutions one can find in \cite{jas, misjas2, KuU2, Mis2, Urbanik93, Urbanik64, Thu, vol2}. \\
 As a simple consequence of Lemma 2.1 and Proposition 2.2 we obtain the following fact:
\begin{prop}\label{prop:3}
Let $\nu \in \mathcal{P}$. For each natural number $n\geqslant 2$ the cumulative distribution function $F_n$ of the measure $\nu^{ \vartriangle_{\alpha} n}$ is equal
$$
F_n(t) = \frac{1}{2\alpha } \left[\alpha \bigl( G(t)^n +1\bigr) + t n  G(t)^{n-1} G'(t) \right], \quad t>0,
$$
and $F(t) = 1 - F(-t)$ for $t<0$, where $G(t) = \widehat{\nu}(1/t)$.
\end{prop}

{\bf Examples.}
\begin{namelist}{ll}
\item{\bf 1.} Let $\nu$ be uniform distribution on $[-1,1]$. Then
\begin{eqnarray*}
G(t) &=& \frac{\alpha}{\alpha+1} |t| \mathbf{1}_{(-1,1)}(t) + \Bigl( |t| - \frac{|t|^{-\alpha}}{\alpha+1} \Bigr) \mathbf{1}_{(1,\infty)}(|t|), \\
F_n(t) &=&  \frac{1}{2}+ \frac{1}{2}\, {\rm sgn}(t) \Bigl(\frac{\alpha}{\alpha+1}\Bigr)^n \left( t + \frac{n}{\alpha} \right)t^{n-1} \pmb{1}_{[0,1)}(t)\\
& & \hspace{-5mm}+ \, \frac{1}{2}\, {\rm sgn}(t) \Bigl(1-\frac{t^{-\alpha}}{\alpha+1}\Bigr)^{n-1}\Bigl(1+\frac{n}{\alpha} + \frac{t^{-\alpha}}{\alpha+1}\left( \frac{n}{t} - 1 \right) \Bigr)\pmb{1}_{[1,\infty)}(t).
\end{eqnarray*}
\item{\bf 2.}
Let $\nu = \widetilde{\delta}_1:= \frac{1}{2}\delta_1 + \frac{1}{2}\delta_{-1}$, thus $G(t) = (1-|t|^{-\alpha})_+$. For each natural number $n\geqslant 2$ the probability measure $\nu^{ \vartriangle_{\alpha} n}$ has the density
$$
dF_n (t) =\frac{\alpha n(n-1)}{2|t|^{2\alpha+1}} \left(1-|t|^{-\alpha}\right)^{n-2} \pmb{1}_{[1,\infty)}(|t|)\, dt.
$$
\item{\bf 3.} Let $\nu = p \widetilde{\delta}_1 + (1-p) \widetilde{\pi}_{p}$, where $p\in(0,1]$ and $\widetilde{\pi}_{p}$ is the symmetric Pareto distribution with density $\widetilde{\pi_{p}}\, (dy) = \alpha |y|^{-p-1} \pmb{1}_{[1,\infty)}(|y|)\, dy$.
    For each natural number $n\geqslant 2$  we have:
\begin{eqnarray*}
\hspace{3mm} F_n(t) & = &  \frac{1}{2}+\frac{1}{2}\,{\rm sgn}(t) \Bigl[1 - \frac{\alpha(1-p)}{\alpha-p}|t|^{-p} + \frac{p(1-\alpha)}{\alpha-p}|t|^{-\alpha}\Bigr]^{n-1} \\
& & \hspace{-9mm} \times \Bigl[ 1+ \frac{(1- p)(np- \alpha)}{\alpha-p}|t|^{-p} -\frac{p(1 -\alpha)(n- 1)}{\alpha-p}|t|^{-\alpha}\Bigr] \pmb{1}_{[1,\infty)}(|t|).
\end{eqnarray*}
\item{\bf 4.}
Let $\nu = \widetilde{\pi}_{2\alpha}$ for $\alpha \in (0,1]$. Since $\widetilde{\delta}_1 \vartriangle_{\alpha} \widetilde{\delta}_1 = \widetilde{\pi}_{2\alpha}$ then using Example 2 we arrive at:
$$
dF_n (t) =\frac{\alpha n(2n-1)}{|t|^{2\alpha+1}} \left(1-|t|^{-\alpha}\right)^{2(n-1)} \pmb{1}_{[1,\infty)}(|t|)\, dt.
$$
\end{namelist}
\section{Kendall random walk}
The direct construction of the symmetric Kendall random walk $\{X_n: n\in\mathbb{N}_0\}$ based on the sequence $(Y_k)$ of i.i.d. unit steps with distribution $\nu$ we have already presented in Definition 1.1. Slightly different direct construction  was given in \cite{KendallWalk}. We recall here one more definition of the Kendall random walk following \cite{BJMR}, where only multidimensional distributions of this process are considered. This approach is more convenient in studying positive and negative excursions, which we consider in this section. 

\begin{definition}
The Kendall random walk or a discrete time L\'{e}vy process under the Kendall generalized distribution is the Markov process $\{X_n: n\in\mathbb{N}_0\}$ with $X_0\equiv 0$ and the transition probabilities 
$$
P_n(x, A) = \mathbf{P}\left\{ X_{n+k} \in A
\big| X_k = x \right\} = \delta_x \vartriangle_{\alpha} \nu^{\vartriangle_{\alpha} n}, \quad n,k \in \mathbb{N},
$$
where measure $\nu \in \mathcal{P}_s$ is called the step distribution.
\end{definition}
 
For any real $a$ we introduce the first hitting times of the half lines $(a,\infty)$ and $(-\infty, a)$ for random walk $\{X_n: n\in\mathbb{N}_0\}$:
$$
\tau_a^+=\min\{n\geqslant 1: X_n>a\}, \quad \tau_a^-=\min\{n\geqslant 1: X_n<a\}
$$
with the convention $\min\emptyset = \infty$. We want to find the joint distribution of the random vector $(\tau_0^+, X_{\tau_0^+})$. In order to attain this we need some lemmas: 

\begin{lem}
\begin{eqnarray*}
h(x,y,t) & := & \delta_x \vartriangle_{\alpha} \delta_y (0,t) = \frac{1}{2} \left( 1 - \left| \frac{xy}{t^2}\right|^{\alpha} \right) \mathbf{1}_{\{|x| < t, |y| < t\}} \\
 & & \hspace{-15mm} = \frac{1}{2}\left[ \Psi\left(\frac{x}{t}\right) + \Psi\left(\frac{y}{t}\right) - \Psi\left(\frac{x}{t}\right) \Psi\left(\frac{y}{t}\right)\right] \mathbf{1}_{\{|x| < t, |y| < t\}}, \\
 \delta_x \vartriangle_{\alpha} \nu \,(0,t) & = & P_1 (x, [0,t)) \\
 & & \hspace{-15mm} = \left[\Psi \left( \frac{x}{t} \right) \Bigl( F(t) - \frac{1}{2} \Bigr) + \frac{1}{2} G(t) - \frac{1}{2} \Psi \left( \frac{x}{t} \right) G(t)\right] \mathbf{1}_{\{|x| < t \}}.
\end{eqnarray*}
\end{lem}

\noindent
{\bf Proof.} With the notation $a = \min\{|x|, |y|\}$, $b= \max\{|x|, |y|\}$ and $\varrho = ({a/b})^{\alpha}$ we have
\begin{eqnarray*}
\lefteqn{ \delta_x \vartriangle_{\alpha} \delta_y (0,t) = T_b \left( (1-\varrho) \widetilde{\delta}_1 + \varrho \widetilde{\pi}_{2\alpha} \right)(0,t)} \\
& = & \frac{1}{2} (1-\varrho)\mathbf{1}(b <t) + \frac{1}{2} \varrho T_b \pi_{2\alpha} (0,t) \\
& = & \frac{1}{2} (1-\varrho)\mathbf{1}(b <t) + \frac{1}{2}\, \varrho\, \mathbf{1}(b <t) \int_1^{t/b} \frac{2\alpha}{s^{2\alpha+1}} ds \\
& = & \frac{1}{2} \Bigl( 1 - \varrho\, \Bigl(\frac{b}{t}\Bigr)^{2\alpha} \Bigr) \mathbf{1}(b <t) = \frac{1}{2} \Bigl( 1 -  \Bigl(\frac{ab}{t^2}\Bigr)^{\alpha} \Bigr) \mathbf{1}(b <t) \\
& = & \frac{1}{2} \Bigl( 1 -  \Bigl|\frac{xy}{t^2}\Bigr|^{\alpha} \Bigr) \mathbf{1}(|x| <t, |y|<t).
\end{eqnarray*}
The final form trivially follows from the identity $1 - ab = 1-a +1-b -(1-a)(1-b)$. The second formula follows from
$$
\delta_x \vartriangle_{\alpha} \nu (0,t) = \int_{\mathbb{R}} \delta_x \vartriangle_{\alpha} \delta_y (0,t) \nu(dy). \eqno{\clubsuit}
$$

The next lemma is crucial for further considerations.

\begin{lem} Let $\{ X_n \colon n \in \mathbb{N}_0 \}$ be the Kendall random walk. Then
\begin{eqnarray*}
\Phi_n(t) & := &  \mathbf{P} \bigl\{ X_1 \leqslant 0, \dots X_{n-1} \leqslant 0, 0< X_n < t \bigr\} \\
 & = & \frac{1}{2^n}\, G(t)^{n-1} \Bigl[ 2n \Bigl( F(t) - \frac{1}{2} \Bigr) - (n-1) G(t)\Bigr].
 \end{eqnarray*}
\end{lem}

\noindent
{\bf Proof.} For $k=1$ we have $\Phi_1(t) = \mathbf{P}\{ 0< X_1 <t\} = F(t) - \frac{1}{2}$. For $k=2$
\begin{eqnarray*}
\Phi_2 (t) & = & \int_{-\infty}^0 \int_0^t P_1(x_1,
 dx_2) P_1(0, dx_1)= \int_{-\infty}^0 P_1(x_1, (0,t))\, \nu(dx_1) \\
& = &  \frac{1}{2} \int_{-t}^0 \Bigl[ 2 \Psi({{x_1}/{t}})\Bigl( F(t) - \frac{1}{2} \Bigr)  + G(t)  - \Psi({{x_1}/{t}}) G(t) \Bigr] \nu(dx_1) \\
 & = & \frac{1}{4}\, \Bigl[ 4G(t)\Bigl( F(t) - \frac{1}{2} \Bigr)  -  G(t)^2 \Bigr].
\end{eqnarray*}
In order to calculate $\Phi_3$ notice first that by Proposition \ref{prop:1} 
\begin{eqnarray*}
\int_{\mathbb{R}} \Psi({{y}/{t}})\left( \delta_{x} \vartriangle_{\alpha} \nu \right)(dy) = \Psi({{x}/{t}}) G(t)
\end{eqnarray*}
and applying the second formula from Lemma 3.1 to the measure $P_1(x_1, \cdot) = \delta_{x_1} \vartriangle_{\alpha} \nu$ we have
\begin{eqnarray*}
\lefteqn{\int_{-\infty}^0 \int_0^t P_1(x_2,  dx_3) P_1(x_1, dx_2) = \int_{-\infty}^0 \left(\delta_{x_2} \vartriangle_{\alpha} \nu\right)(0,t) \, \left(\delta_{x_1} \vartriangle_{\alpha} \nu\right) (dx_2)} \\
&& \hspace{-5mm} = \frac{1}{2}\int_{-t}^0 \left[ 2 \Psi({{x_2}/{t}})\Bigl(\! F(t) - \frac{1}{2} \Bigr)  + G(t)  - \Psi({{x_2}/{t}}) G(t) \right] \left(\delta_{x_1} \vartriangle_{\alpha} \nu\right) (dx_2) \\
&& \hspace{-5.5mm}= \frac{1}{2}\left[ \Psi({{x_1}/{t}}) G(t)\! \Bigl(\! F(t) - \frac{1}{2} \Bigr)+ G(t)\delta_{x_1} \vartriangle_{\alpha} \nu(0,t)  - \frac{1}{2}\Psi({{x_1}/{t}}) G(t)^2 \right] \\
&& \hspace{-5mm}= \frac{1}{4}\left[ 4 \Psi({{x_1}/{t}}) G(t)\! \Bigl(\! F(t) - \frac{1}{2} \Bigr)+ G(t)^2  - 2 \Psi({{x_1}/{t}}) G(t)^2 \right].
\end{eqnarray*}
Consequently
\begin{eqnarray*}
\lefteqn{\Phi_3 (t)= \int_{-\infty}^0 \int_{-\infty}^0 \int_0^t P_1(x_2,  dx_3) P_1(x_1, dx_2) P_1(0,dx_1) } \\
 && \hspace{-10mm} = \frac{1}{4}\int_{-t}^0 \left[ 4 \Psi({{x_1}/{t}}) G(t)\! \Bigl(\! F(t) - \frac{1}{2} \Bigr) + G(t)^2  - 2 \Psi({{x_1}/{t}}) G(t)^2 \right] \nu(dx_1) \\
 && \hspace{-10mm} = \frac{1}{4} \left[ 2 G(t)^2 \! \Bigl(\! F(t) - \frac{1}{2} \Bigr)+ G(t)^2 \! \Bigl(\! F(t) - \frac{1}{2} \Bigr) - G(t)^3 \right] \\
 & = & \frac{1}{2^3} \left[ 6 G(t)^2 \! \Bigl(\! F(t) - \frac{1}{2} \Bigr) - 2 G(t)^3 \right].
\end{eqnarray*}
Simple, but laborious, application of mathematical induction ends the proof. \qed

\begin{lem} The random variable $\tau_0^+$ (and, by symmetry of the Kendall random walk, also variable $\tau_0^-$) has geometric distribution $P(\tau_0^{+} = k) = \frac{1}{2^k}, \quad k=1,2,\cdots$  it has the following moments generating function:
$$
\mathbf{E} s^{\tau_0^+} = \frac{\frac{s}{2}}{ 1 - \frac{s}{2}}, \quad 0\leqslant s<2.
$$
\end{lem}

\noindent
{\bf Proof.}
Notice that $\lim_{t\rightarrow \infty} G(t) = 1$ since $\lim_{t\rightarrow \infty} \Psi(1/t) = 1$. Now it is enough to apply Lemma 3.1:
$$
P(\tau_0^{+} = k)  =   \mathbf{P} \left\{  X_1 \leqslant 0, \cdots, X_{k-1} \leqslant 0, X_k > 0 \right\}  = \Phi_k (\infty) = \frac{1}{2^k}. \eqno{\clubsuit}
$$
\begin{lem}
$$
\mathbf{P}\left\{ X_{\tau_0^{+}} < t \right\} = \frac{ 4F(t) - 2 - G(t)^2}{(2-G(t))^2}.
$$
\end{lem}

\noindent
{\bf Proof.} It is enough to apply Lemma 3.1:
\begin{eqnarray}
\lefteqn{ \mathbf{P}\left\{ X_{\tau_0^{+}} < t \right\}   =
\sum\limits_{k=1}^{\infty} \mathbf{P} \left\{  X_{\tau_0^{+}} < t,  \tau_0^{+} =k\right\}= \sum\limits_{k=1}^{\infty} \Phi_k (t)} \notag \\
& = & \sum_{n=1}^{\infty} \frac{1}{2^n}\, G(t)^{n-1} \Bigl[ 2n \Bigl( F(t) - \frac{1}{2} \Bigr) - (n-1) G(t)\Bigr] \notag \\
& = & \frac{1}{(2-G(t))^2}\, \Bigl[ 4F(t) - 2 - G(t)^2 \Bigr]. \hspace{40mm}{\clubsuit}  \notag
 \end{eqnarray}

\begin{definition}
Let $\{X_n \colon n \in \mathbb{N}_0\}$ be the Kendall random walk independent of the random variable  $ N_{s/2}$ with the geometric distribution $\mathbf{P} \{  N_{s/2} = k\} = (\frac{s}{2})^{k-1} (1 - \frac{s}{2})$, $0\leqslant s \leqslant 1$, for $k=1,2,\dots$. By the geometric-Kendall random variable we understand the following
$$
Z_{s/2} = \sum_{k=1}^{\infty} X_k \mathbf{1}_{\{ N_{s/2} = k\}}.
$$
\end{definition}
Notice that $\mathbf{E} \Psi\left( {{X_k}/{t}}\right) = G(t)^k$ by Proposition 2.2, thus we have
\begin{eqnarray*}
\mathbf{E} \Psi({{Z_{s/2}}/t}) & = & \sum\limits_{k=1}^{\infty} \mathbf{E} \Psi\left( {{X_k}/{t}}\right) \left(\frac{s}{2}\right)^{k-1}\Bigl( 1 - \frac{s}{2} \Bigr)   \\
& = &  \sum\limits_{k=1}^{\infty} G(t)^k \left( \frac{s}{2} \right)^{k-1} \!\left( 1 - \frac{s}{2} \right)
= \frac{G(t) \bigl( 1- \frac{s}{2} \bigr)}{1 - \frac{s}{2}\, G(t)}.
\end{eqnarray*}
We are ready now to prove the main result of the paper: the Wiener-Hopf factorization for the Kendall random walk. This gives the Laplace-Williamson transform $\mathbf{E} s^{\tau_0^+} \Psi ( uX_{\tau_0^+})$, thus consequently, also the joint distribution of $({\tau_0^+}, X_{\tau_0^+})$.
\begin{thm}
Let $\{X_n: n\in\mathbb{N}_0\}$ be a random walk under the Kendall convolution with unit step $X_1 \sim F$ such that $F(0)=\frac{1}{2}$. Then
$$
\mathbf{E} s^{\tau_0^+}\Psi\left(uX_{\tau_0^+}\right) = \frac{\frac{s}{2}}{1-\frac{s}{2} } \cdot \frac{(1-\frac{s}{2}) G({1/u})}{1-\frac{s}{2} G({1/u})} = \mathbf{E} s^{\tau_0^+}\, \mathbf{E}\Psi\left(uZ_{s/2}\right).
$$
\end{thm}
\noindent {\bf Proof.} Let $H(s,u) := \mathbf{E} s^{\tau_0^+} \Psi \left( uX_{\tau_0^+} \right)$. Then
\begin{eqnarray*}
H(s,u) & = & \sum_{k=1}^{\infty}  s^{k}\int_{0}^{\infty} \Psi(ux) d\Phi_k(x) = \int\limits_{0}^{\infty} \! \Psi(ux)\, d \biggl( \sum_{k=1}^{\infty}  s^{k} \Phi_k(x)\biggr)_x \\
&& \hspace{-15mm} =
\int \limits_{0}^{\infty} \!\!\Psi(ux)\, d \biggl( \sum_{k=1}^{\infty} \Bigl[\frac{s G(x)}{2}\Bigr]^{k} \!\! \frac{1}{G(x)} \left[ 2k \Bigl[ F(x) - \frac{1}{2} \Bigr]\! - (k-1) G(x)\right]\biggr)_x \\
&& \hspace{-15mm} = 
\int_{0}^{\infty}\!\!\Psi(ux)\, d \biggl[ \frac{s (F(x) - \frac{1}{2}) - \frac{1}{4} s^2 G(x)^2}{\left(1-\frac{1}{2} s G(x)\right)^2 }\biggr]_x \\
&& \hspace{-15mm} = \alpha u^{\alpha} \int_{0}^{1/u}\!\!x^{\alpha - 1} \, \frac{s (F(x) - \frac{1}{2}) - \frac{1}{4} s^2 G(x)^2}{\left(1-\frac{1}{2} s G(x)\right)^2 }\,\, dx.
\end{eqnarray*}
In these calculations we write $d(f(s,x))_x$ for the differential of the function $f(s,x)$ with respect to $x$. To get the last expression we need to know that $G(0) = 0$.
Now we calculate $G'(x)$. Since
$$
G(x) = \int_{\mathbb{R}} \Psi \left( \frac{z}{x}\right) \nu(dz) = 2 \int_0^{x} \left( 1 - \Bigl|\frac{z}{x}\Bigr|^{\alpha} \right) \nu(dz),
$$
then integrating by parts we obtain
$$
\left( G(x) + 1 \right) x^{\alpha} = 2 \alpha \int_0^x z^{\alpha - 1} F(z) dz.
$$
Now it is enough to differentiate both sides of this equality to obtain
$$
G'(x) = 2 \alpha x^{-1} \left( F(x) - \frac{1}{2} - \frac{1}{2} G(x) \right).
$$
We see that
$$
H(s,u) = u^{\alpha} \int_{0}^{1/u}\left[  \frac{x^{\alpha} \frac{s}{2} G(x)}{1-\frac{s}{2} G(x)}\right]' dx =  \frac{\frac{s}{2} G({1/u})}{1-\frac{s}{2} G({1/u})}. \eqno{\clubsuit}
$$
The same method as in the proof of the last theorem we are using in the next lemma.

\begin{lem}
$$
\mathbf{E} X_{\tau_0^{+}}^{\alpha} = 2 \mathbf{E} | Y_1 |^{\alpha}
$$
\end{lem}

\noindent {\bf Proof.} In the following calculations we take $R>0$:
\begin{eqnarray*}\lefteqn{\mathbf{E} X_{\tau_0^{+}}^{\alpha} = \int_0^{\infty} x^{\alpha} d \left( \frac{ F(x) - \frac{1}{2} - \frac{1}{4} G(x)^2}{(1 - \frac{1}{2} G(x))^2}\right)_x } \\
 & = & \lim_{R\rightarrow \infty} \left[ x^{\alpha} \frac{ F(x) - \frac{1}{2} - \frac{1}{4} G(x)^2}{(1 - \frac{1}{2} G(x))^2} \Big|_0^R - \int_0^R \Bigl[ \frac{ x^{\alpha} \frac{1}{2} G(x)}{(1 - \frac{1}{2} G(x))}\Bigr]' dx \right] \\
 & = & \lim_{R\rightarrow \infty} \left[ x^{\alpha} \frac{ F(x) - \frac{1}{2} - \frac{1}{2} G(x)}{(1 - \frac{1}{2} G(x))^2} \Big|_0^R  \right] \\
 & = & \lim_{R\rightarrow \infty} \left[ \frac{ x^{\alpha}}{(1 - \frac{1}{2} G(x))^2} \Bigl[ \int_0^x \nu(du) - \int_0^x \Bigl[ 1 - \left(\frac{u}{x}\right)^{\alpha}\Bigr] \nu(du)\Bigr]  \right]_0^R \\
 & = & \lim_{R\rightarrow \infty} \left[ \frac{ 1}{(1 - \frac{1}{2} G(x))^2} \int_0^x u^{\alpha}  \nu(du) \right]_0^R = 4 \int_0^{\infty} u^{\alpha}  \nu(du),
\end{eqnarray*}
which was to be shown. Notice that the assumption $\mathbf{E}|Y_1|^{\alpha} < \infty$ is irrelevant in these calculations.  \qed

\end{document}